\documentclass[12pt]{amsart}

\input epsf
\theoremstyle{plain}
\newtheorem{thm}[subsubsection]{Theorem}
\newtheorem{lem}[subsubsection]{Lemma}
\newtheorem{prop}[subsubsection]{Proposition}
\newtheorem{cor}[subsubsection]{Corollary}

\theoremstyle{definition}
\newtheorem{rem}[subsubsection]{Remark}
\newtheorem{defn}[subsubsection]{Definition}
\newtheorem{ex}[subsubsection]{Example}

\def\w{\wedge}

\def\ui{\underline{i}}

\def\res{|}

\def\L\mathcal{{L}}
\def\A\mathcal{{A}}
\def\a{\alpha}
\def\b{\beta}
\def\g{\gamma}

\def\e{\epsilon}

\def\L{\Lambda}
\def\m{\mu}
\def\n{\nu}
\def\o{\omega}

\def\t{\tau}
\def\te{\theta}

\def\ui{\underline{i}}
\def\uj{\underline{j}}

\def\uk{\underline{k}}
\def\ue{\underline{e}}
\def\uf{\underline{f}}

\def\ni{\noindent}

\newcommand{\AAA}{\mathbb{A}}
\newcommand{\PPP}{\mathbb{P}}

\DeclareMathOperator{\MMod}{\mathrm{mod }}
\DeclareMathOperator{\Mult}{\mathrm{mult}}

\DeclareMathOperator{\Sing  }{\mathrm{Sing}}
\begin{document}

\title{Richardson Varieties in the Grassmannian}
\author {Victor Kreiman}
\address{Department of Mathematics\\ Northeastern University\\ Boston, MA
02115} \email{vkreiman@lynx.neu.edu}
\author[V. Lakshmibai]{V. Lakshmibai${}^{\dag}$}
\address{Department of Mathematics\\ Northeastern University\\ Boston, MA
02115} \email{lakshmibai@neu.edu}
\dedicatory{Dedicated to
Professor J. Shalika on his sixtieth birthday}
\thanks{${}^{\dag}$ Partially supported
by NSF grant DMS 9971295.}

\begin{abstract}
The Richardson variety $X_w^v$ is defined to be the intersection
of the Schubert variety $X_w$ and the opposite Schubert variety
$X^v$. For $X_w^v$ in the Grassmannian, we obtain a standard
monomial basis for the homogeneous coordinate ring of $X_w^v$. We
use this basis first to prove the vanishing of $H^i(X_w^v,L^m)$,
$i > 0 $, $m \geq 0$, where $L$ is the restriction to $X_w^v$ of
the ample generator of the Picard group of the Grassmannian; then
to determine a basis for the tangent space and a criterion for
smoothness for $X_w^v$ at any $T$-fixed point $e_\t$; and finally
to derive a recursive formula for the multiplicity of $X_w^v$ at
any $T$-fixed point $e_\t$. Using the recursive formula, we show
that the multiplicity of $X_w^v$ at $e_\t$ is the product of the
multiplicity of $X_w$ at $e_\t$ and the multiplicity of $X^v$ at
$e_\t$. This result allows us to generalize the
Rosenthal-Zelevinsky determinantal formula for multiplicities at
$T$-fixed points of Schubert varieties to the case of Richardson
varieties.
\end{abstract}
\maketitle

\setcounter{tocdepth}{1} \tableofcontents

\section*{Introduction}

Let $G$ denote a semisimple, simply connected, algebraic group
defined over an algebraically closed field $K$ of arbitrary
characteristic. Let us fix a maximal torus $T$ and a Borel
subgroup $B$  containing $T$. Let $W$ be the Weyl group ($
N(T)/T$,  $N(T)$ being the normalizer of $T$). Let $Q$ be a
parabolic subgroup  of $G$ containing $B$, and $W_Q$, the Weyl
group of $Q$. For the action of $G$ on $G/Q$ given by left
multiplication, the $T$-fixed points are precisely the cosets
$e_w:=wQ$ in $G/Q$. For $w \in W/W_Q$, let $X_w$ denote the {\it
Schubert variety} (the Zariski closure of the $B$-orbit $Be_w$ in
$G/Q$ through the $T$-fixed point $e_w$), endowed with the
canonical structure of a closed, reduced subscheme of $G/Q$. Let
$B^-$ denote the Borel subgroup of $G$ opposite to $B$ (it is the
unique Borel subgroup of $G$ with the property $B\cap B^-=T$). For
$v \in W/W_Q$, let $X^v$ denote the {\it opposite Schubert
variety}, the Zariski closure of the $B^-$-orbit $B^-e_v$ in
$G/Q$.

Schubert and opposite  Schubert varieties play an important role
in the study of the generalized flag variety $G/Q$, especially,
the algebraic-geometric and representation-theoretic aspects of
$G/Q$. A more general class of subvarieties in $G/Q$ is the class
of Richardson varieties; these are varieties of the form $X^v_w:=
X_w\cap X^v $, the intersection of the Schubert variety $X_w$ with
opposite Schubert variety $X^v$. Such varieties were first
considered by Richardson in (cf. \cite{Ri}), who shows that such
intersections are reduced and irreducible. Recently, Richardson
varieties have shown up in several contexts: such double coset
intersections $BwB\cap B^-xB$ first appear in \cite{K-Lu$_1$},
\cite{K-Lu$_2$}, \cite{Ri}, \cite{Ri-Ro-S}. Very recently,
Richardson varieties have also appeared in the context of K-theory
of flag varieties (\cite{Br}, \cite{L-L}). They also show up in
the construction of certain degenerations of Schubert varieties
(cf. \cite{Br}).

In this paper, we present results for Richardson varieties in the
Grassmannian variety. Let $G_{d,n}$ be the Grassmannian variety of
$d$-dimensional subspaces of $K^n$, and $p:G_{d,n}\hookrightarrow
\mathbb{P}^N \,(=\mathbb{P}(\wedge^dK^n))$, the Pl\"ucker
embedding (note that $G_{d,n}$ may be identified with $G/P,
G=SL_n(K), P$ a suitable maximal parabolic subgroup of $G$). Let
$X:=X_w\cap X^v$ be a Richardson variety in $G_{d,n}$. We first
present a Standard monomial theory for $X$ (cf. Theorem
\ref{mainth}). Standard monomial theory (SMT) consists in
constructing an explicit basis for the homogeneous coordinate ring
of $X$. SMT for Schubert varieties was first developed by the
second author together with Musili and Seshadri in a series of
papers, culminating in \cite{g/p5}, where it is established for
all classical groups. Further results concerning certain
exceptional and Kac--Moody groups led to conjectural formulations
of a general SMT, see \cite{ls2}. These conjectures were then
proved by Littelmann, who introduced new combinatorial and
algebraic tools: the path model of representations of any
Kac-Moody group, and Lusztig's Frobenius map for quantum groups at
roots of unity (see \cite{Li$_1$,Li$_2$}); recently, in
collaboration with Littelmann (cf. \cite{L-L}), the second author
has extended the results of \cite{Li$_2$} to Richardson varieties
in $G/B$, for any semisimple $G$. Further, in collaboration with
Brion (cf. \cite{b-l}), the second author has also given a purely
geometric construction of standard monomial basis for Richardson
varieties in $G/B$, for any semisimple $G$; this construction in
loc. cit. is done using certain flat family with generic fiber
$\cong$ diag$(X^v_w)\subset X^v_w\times X^v_w$, and the special
fiber $\cong \cup_{v\le x\le w}X^v_x\times X_w^x$.

If one is concerned with just Richardson varieties in the
Grassmannian, one could develop a SMT in the same spirit as in
\cite{Mu} using just the Pl\"ucker coordinates, and one doesn't
need to use any quantum group theory nor does one need the
technicalities of \cite{b-l}. Thus we give a self-contained
presentation of SMT for unions of Richardson varieties in the
Grassmannian. We should remark that Richardson vatieties in the
Grassmannian are also studied in \cite{st}, where these varieties
are called {\it skew Schubert varieties}, and standard monomial
bases for these varieties also appear in loc. cit. (Some
discussion of these varieties also appears in \cite{hp2}.) As a
consequence of our results for unions of Richardson varieties, we
deduce the vanishing of $H^i(X,L^m), i\ge 1, m\ge 0, L$, being the
restriction to $X$ of ${\mathcal O}_{\mathbb{P}^N}(1)$ (cf.
Theorem \ref{main}); again, this result may be deduced using the
theory of Frobenius-splitting (cf. \cite{m-r}), while our approach
uses just the classical Pieri formula. Using the standard monomial
basis, we then determine the tangent space and also the
multiplicity at any $T$-fixed point $e_\tau$ on $X$. We first give
a recursive formula for the multiplicity of $X$ at $e_\tau$ (cf.
Theorem \ref{8.8.5.3}). Using the recursive formula, we derive a
formula for the multiplicity of $X$ at $e_\tau$ as being the
product of the multiplicities at $e_\tau$ of $X_w$ and $X^v$ (as
above, $X=X_w\cap X^v$) (cf. Theorem \ref{prodmult}). Using the
product formula, we get a generalization of Rosenthal-Zelevinsky
determinantal formula (cf. \cite{rose-zel}) for the multiplicities
at singular points of Schubert varieties to the case of Richardson
varieties (cf. Theorem \ref{rzr}).   It should be mentioned that
the multiplicities of Schubert varieties at $T$-fixed points
determine their multiplicities at all other points, because of the
$B$-action; but this does not extend to Richardson varieties,
since Richardson varieties have only a $T$-action. Thus even
though, certain smoothness criteria at $T$-fixed points on a
Richardson variety are given in Corollaries \ref{smcrit1} and
\ref{smcrit2}, the problem of the determination of singular loci
of Richardson varieties still remains open.

In \S \ref{gdn}, we present basic generalities on the Grassmannian
variety and the Pl\"ucker embedding.  In \S \ref{sor}, we define
Schubert varieties, opposite Schubert varieties, and the more
general Richardson varieties in the Grassmannian and give some of
their basic properties. We then develop a standard monomial theory
for a Richardson variety $X_w^v$ in the Grassmannian in \S
\ref{smt1} and extend this to a standard monomial theory for
unions and nonempty intersections of Richardson varieties in the
Grassmannian in \S \ref{unn}.  Using the standard monomial theory,
we obtain our main results in the three subsequent sections.  In
\S \ref{van}, we prove the vanishing of $H^i(X_w^v,L^m)$, $i > 0
$, $m \geq 0$, where $L$ is the restriction to $X_w^v$ of the
ample generator of the Picard group of the Grassmannian.  In \S
\ref{tss}, we determine a basis for the tangent space and a
criterion for smoothness for $X_w^v$ at any $T$-fixed point
$e_\t$. Finally, in \S \ref{mas}, we derive several formulas for
the multiplicity of $X_w^v$ at any $T$-fixed point $e_\t$.

We are thankful to the referee for many valuable comments and
suggestions, especially for the alternate proof of Theorem
\ref{prodmult}.

\section{The Grassmannian Variety $G_{d, n}$}\label{gdn}

Let $K$ be the base field, which we assume to be algebraically
closed of arbitrary characteristic.  Let $d$ be such that $1 \le d
< n$. The {\it Grassmannian} $G_{d,n}$ is the set of all
$d$-dimensional subspaces of $K^n$.  Let $U$ be an element of
$G_{d,n}$ and $\{ a_1, \ldots a_d \}$ a basis of $U$, where each
$a_j$ is a vector of the form \vspace{-.05in}
\begin{displaymath}
a_j = \left(\begin{array}{c} a_{1j} \\ a_{2j} \\ \vdots \\ a_{nj}
\end{array}\right),\, \textrm{with } a_{ij} \in K.
\end{displaymath}
Thus, the basis $\{a_1, \cdots, a_d\}$ gives rise to an $n \times
d$ matrix $A = (a_{ij})$ of rank $d$, whose columns are the
vectors $a_1, \cdots, a_d$.

We have a canonical embedding
\begin{equation*}
p:G_{d,n}\hookrightarrow \mathbb{P}(\wedge^dK^n) \ , \ U \mapsto
[a_1 \wedge \cdots \wedge a_d]
\end{equation*}
called the {\it Pl\" ucker embedding}. It is well known that $p$
is a closed immersion; thus $G_{d, n}$ acquires the structure of a
projective variety. Let
\begin{equation*}
I_{d,n}=\{\underline{i}=(i_1,\dots,i_d) \in \mathbb{N}^d \, : \,
1\le i_1<\dots<i_d\le n\}.
\end{equation*}
Then the projective coordinates ({\it Pl\"ucker coordinates}) of
points in $\mathbb{P}(\wedge^d K^n)$ may be indexed by $I_{d,n}$;
for $\underline{i} \in I_{d,n}$, we shall denote the
$\underline{i}$-th component of $p$ by $p_{\underline{i}}$, or
$p_{i_1, \cdots, i_d}$.  If a point $U$ in $G_{d, n}$ is
represented by the $n \times d$ matrix $A$ as above, then $p_{i_1,
\cdots, i_d}(U) = \textrm{det}(A_{i_1, \ldots, i_d})$, where
$A_{i_1, \ldots, i_d}$ denotes the $d \times d$ submatrix whose
rows are the rows of $A$ with indices $i_1, \ldots, i_d$, in this
order.

For $\ui\in I_{d,n}$ consider the point $e_{\ui}$ of $G_{d,n}$
represented by the $n\times d$ matrix whose entries are all 0,
except the ones in the $i_j$-th row and $j$-th column, for each
$1\le j \le d$, which are equal to 1. Clearly, for $\ui, \uj\in
I_{d,n}$,
\begin{equation*}
p_{\ui}(e_{\uj})=
\begin{cases}
1,&\text{\rm if $\ui=\uj$;}\\ 0,&\text{\rm otherwise.}
\end{cases}
\end{equation*}

We define a partial order $\ge$ on $I_{d,n}$ in the following
manner: if $\ui = (i_1, \cdots, i_d)$ and $\uj = (j_1, \cdots,
j_d)$, then  $\ui\ge \uj \Leftrightarrow i_t\ge j_t , \forall t$.
The following well known theorem gives the defining relations of
$G_{d, n}$ as a closed subvariety of $\mathbb{P}(\wedge^dK^n)$
(cf. \cite{hodge}; see \cite{L-G} for details):

\begin{thm}\label{T:greq}
The Grassmannian $G_{d,n}\subset \mathbb{P}(\wedge^dK^n)$ consists
of the zeroes in

\ni $\mathbb{P}(\wedge^dK^n)$ of quadratic polynomials of the form
\begin{equation*}
p_{\ui}p_{\uj}-\sum \pm p_\a p_\b
\end{equation*}
for all $\ui ,\, \uj \in I_{d, n}$, $\ui,\, \uj$ non-comparable,
where $\a, \, \b$ run over a certain subset of $I_{d, n}$ such
that $\a
>$ both $\ui$ and $\uj$, and $\b <$ both $\ui$ and $\uj$.
\end{thm}

\subsection{Identification of $G/P_d$ with $G_{d,n}$}\label{fd}
Let $G=SL_n(K)$.  Let $P_d$ be the maximal parabolic subgroup
\begin{equation*}
P_d=\left\{ A\in G\biggm| A=
\begin{pmatrix}
*&*\\ 0_{(n-d)\times d}&*
\end{pmatrix}\right\}.
\end{equation*}
For the natural action of $G$ on $\mathbb{P}(\wedge^d K^n)$, we
have, the isotropy at $[e_1\wedge \cdots \wedge e_d]$ is $P_d$
while the orbit through $[e_1\wedge \cdots\wedge e_d]$ is
$G_{d,n}$. Thus we obtain a surjective morphism $\pi:G\rightarrow
G_{d,n},\ g\mapsto g\cdot a$, where $a=[e_1\wedge \cdots \wedge
e_d]$. Further, the differential $d\pi_e :$ Lie$G \rightarrow
T(G_{d,n})_a$ (= the tangent space to $G_{d,n}$ at $a$) is easily
seen to be surjective. Hence we obtain an identification
$f_d:G/P_d\cong G_{d,n}$ (cf. \cite{borel1}, Proposition 6.7).

\subsection{Weyl Group and Root System}\label{weylroot} Let $G$
and $P_d$ be as above. Let $T$ be the subgroup of diagonal
matrices in $G$, $B$ the subgroup of upper triangular matrices in
$G$, and $B^-$ the subgroup of lower triangular matrices in $G$.
Let $W$ be the Weyl group of $G$ relative to $T$, and $W_{P_d}$
the Weyl group of $P_d$. Note that $W = S_n$, the group of
permutations of a set of $n$ elements, and that $W_{P_d} = S_d
\times S_{n-d}$. For a permutation $w$ in $S_n$, $l(w)$ will
denote the usual length function. Note also that $I_{d, n}$ can be
identified with $W / W_{P_d}$. In the sequel, we shall identify
$I_{d,n}$ with the set of ``minimal representatives" of $W /
W_{P_d}$ in $S_n$; to be very precise, a $d$-tuple $\ui\in
I_{d,n}$ will be identified with the element
$(i_1,\dots,i_d,j_1,\dots,j_{n-d})\in\mathcal{S}_n$, where
$\{j_1,\dots,j_{n-d}\}$ is the complement of $\{i_1,\dots,i_d\}$
in $\{1,\dots,n\}$ arranged in increasing order.  We denote the
set of such minimal representatives of $S_n$ by $W^{P_d}$.

Let $R$ denote the root system of $G$ relative to $T$, and $R^+$
the set of positive roots relative to $B$.  Let $R_{P_d}$ denote
the root system of $P_d$, and $R_{P_d}^+$ the set of positive
roots.

\section{Schubert, Opposite Schubert, and Richardson Varieties in
$G_{d,n}$}\label{sor} For $1\le t\le n$, let $V_t$ be the subspace
of $K^n$ spanned by $\{ e_1,\dots,e_t\}$, and let $V^t$ be the
subspace spanned by $\{ e_n, \ldots, e_{n-t+1}\}$. For each
$\ui\in I_{d,n}$, the {\em Schubert variety $X_{\ui}$} and {\em
Opposite Schubert variety $X^{\ui}$} associated to $\ui$ are
defined to be
\[
X_{\ui}=\{U\in G_{d,n}\mid \dim (U\cap V_{i_t})\ge t\ ,\ 1\le t\le
d\},
\]
\[
X^{\ui}=\{U\in G_{d,n}\mid \dim (U\cap V^{n-i_{(d-t+1)}+1})\ge t\
,\ 1\le t\le d\}.
\]
For ${\ui}, {\uj} \in I_{d, n}$, the {\em Richardson Variety
$X_{\ui}^{\uj}$}is defined to be $X_{\ui} \cap X^{\uj}$.  For
${\ui}, {\uj}, {\ue}, {\uf} \in I_{d,n}$, where ${\ue}=(1, \ldots,
d)$ and ${\uf}=(n+1-d, \ldots, n)$, note that
$G_{d,n}=X_{\uf}^{\ue}$, $X_{\ui}=X_{\ui}^{\ue}$, and
$X^{\ui}=X_{\uf}^{\ui}$.

For the action of $G$ on $\mathbb{P}(\wedge^dK^n)$,  the $T$-fixed
points are precisely the points corresponding to the
$T$-eigenvectors in $\wedge^dK^n$. Now
\[
\wedge^dK^n=\bigoplus_{\ui\in I_{d,n}} Ke_{\ui} \text{\rm ,\quad
as $T$-modules,}
\]
where for $\ui =(i_1,\cdots ,i_d),\ e_{\ui}=e_{i_1}\wedge \cdots
\wedge e_{i_d}$. Thus the $T$-fixed points in
$\mathbb{P}(\wedge^dK^n)$ are precisely $\left[ e_{\ui}\right]$,
$\ui\in I_{d,n}$, and these points, obviously, belong to
$G_{d,n}$. Further, the Schubert variety $X_{\ui}$ associated to
$\ui$ is simply the Zariski closure of the $B$-orbit $B\left[
e_{\ui}\right]$ through the $T$-fixed point $\left[
e_{\ui}\right]$ (with the canonical reduced structure), $B$ being
as in \S \ref{weylroot}. The opposite Schubert variety $X^{\ui}$
is the Zariski closure of the $B^-$-orbit $B^-\left[
e_{\ui}\right]$ through the $T$-fixed point $\left[
e_{\ui}\right]$ (with the canonical reduced structure), $B^-$
being as in \S \ref{weylroot}.

\subsection{Bruhat Decomposition}\label{bruh} Let $V=K^n$.  Let $\ui\in I_{d,n}$. Let
$C_{\ui}=B[e_{\ui}]$ be the {\em Schubert cell} and
$C^{\ui}=B^-[e_{\ui}]$ the {\em opposite Schubert cell} associated
to $\ui$. The $C_{\ui}$'s provide a cell decomposition of
$G_{d,n}$, as do the $C^{\ui}$'s. Let $X=V\oplus\dots\oplus V$ (d
times). Let

$$\pi:X\to\w^dV, (u_1,\dots,u_d)\mapsto u_1\w\dots\w u_d,$$ and
$$p:\w^dV\setminus\{0\}\to\mathbb{P}(\w^dV), u_1\w\dots\w
u_d\mapsto[u_1\w\dots\w u_d].$$ Let $v_{\ui}$ denote the point
$(e_{i_1},\dots,e_{i_d})\in X$.

Identifying $X$ with $M_{n\times d}$, $v_{\ui}$ gets identified
with the $n\times d$ matrix whose entries are all zero except the
ones in the $i_j$-th row and $j$-th column, $1\le j\le d$, which
are equal to $1$. We have

$$B\cdot v_{\ui}=\{ A\in M_{n\times d}\mid x_{ij}=0,i>i_j,\textrm{
and } \prod_t x_{i_t t}\ne 0\},$$

$$B^-\cdot v_{\ui}=\{ A\in M_{n\times d}\mid
x_{ij}=0,i<i_j,\textrm{ and } \prod_t x_{i_t t}\ne 0\}.$$ Denoting
$\overline{B\cdot v_{\ui}}$ by $D_{\ui}$, we have $D_{\ui}=\{A\in
M_{n\times d}\mid x_{ij}=0,i>i_j\}$. Further, $\pi(B\cdot
v_{\ui})=p^{-1}(C_{\ui})$, $\pi(D_{\ui})=\widehat{X_{\ui}}$, the
cone over $X_{\ui}$. Denoting $\overline{B^-\cdot v_{\ui}}$ by
$D^{\ui}$, we have $D^{\ui}=\{A\in M_{n\times d}\mid
x_{ij}=0,i<i_j\}$. Further, $\pi(B^-\cdot
v_{\ui})=p^{-1}(C^{\ui})$, $\pi(D^{\ui})=\widehat{X^{\ui}}$, the
cone over $X^{\ui}$. From this, we obtain

\begin{thm}\label{less}
\begin{enumerate}
\item Bruhat Decomposition:
$X_{\uj}=\dot{\bigcup\limits_{\ui\,\le\,\uj}}Be_{\ui}$,
$X^{\uj}=\dot{\bigcup\limits_{\ui\,\ge\,\uj}}B^-e_{\ui}$.
\item $X_{\ui}\subseteq X_{\uj}$ if and only if $\ui\le\uj$.
\item $X^{\ui}\subseteq X^{\uj}$ if and only if $\ui\ge\uj$.
\end{enumerate}
\end{thm}

\begin{cor}\label{morerich}
\begin{enumerate}
\item $X_{\uj}^{\uk}$ is nonempty $\iff \uj \geq \uk$; further,
when $X_{\uj}^{\uk}$ is

\ni nonempty, it is reduced and irreducible of dimension
$l(w)-l(v)$, where

\ni $w$ (resp. $v$) is the permutation in $S_n$ representing $\uj$
(resp. $\uk$) as in \S \ref{weylroot}.
\item $p_{\uj}\bigr|_{X_{\ui}^{\uk}}\ne
0\iff\ui\ge\uj\ge\uk.$
\end{enumerate}
\end{cor}
\begin{proof}

\ni (1) Follows from \cite{Ri}.  The criterion for $X_{\uj}^{\uk}$
to be nonempty, the irreducibility, and the dimension formula are
also proved in \cite{Deo}.

\ni (2) From Bruhat decomposition, we have
$p_{\uj}\res_{X_{\ui}}\ne 0 \iff e_{\uj}\in X_{\ui}$; we also have
$p_{\uj}\res_{X^{\uk}}\ne 0 \iff e_{\uj}\in X^{\uk}$.  Thus
$p_{\uj}\res_{X_{\ui}^{\uk}}\ne 0 \iff e_{\uj}\in X_{\ui}^{\uk}$.
Again from Bruhat decomposition, we have $e_{\uj}\in X_{\ui}^{\uk}
\iff \ui\ge\uj\ge\uk$. The result follows from this.
\end{proof}
For the remainder of this paper, we will assume that all our
Richardson varieties are nonempty.

\begin{rem}
In view of Theorem \ref{less}, we have $X_{\ui}\subseteq X_{\uj}$
if and only if $\ui\le \uj$. Thus, under the set-theoretic
bijection between the set of Schubert varieties and the set
$I_{d,n}$, the partial order on the set of Schubert varieties
given by inclusion induces the partial order $\ge $ on $I_{d,n}$.
\end{rem}

\subsection{More Results on Richardson Varieties}\label{mrrch}

\begin{lem}\label{bstable}
Let $X \subseteq G_{d, n}$ be closed and $B$-stable (resp.
$B^-$-stable).  Then $X$ is a union of Schubert varieties (resp.
opposite Schubert varieties).
\end{lem}
The proof is obvious.

\begin{lem}\label{inter}
Let $X_1$, $X_2$ be two Richardson varieties in $G_{d,n}$ with
nonempty intersection. Then $X_1\cap X_2$ is a Richardson variety
(set-theoretically).
\end{lem}
\begin{proof}
We first give the proof when $X_1$ and $X_2$ are both Schubert
varieties. Let $X_1=X_{\t_1}$, $X_2=X_{\t_2}$, where
$\t_1=(a_1,\dots, a_d)$, $\t_2=(b_1,\dots, b_d)$. By Lemma
\ref{bstable}, $X_1\cap X_2=\cup X_{w_i}$, where $w_i<\t_1$,
$w_i<\t_2$. Let $c_j=\min\{a_j,b_j\}$, $1\le j\le d$, and
$\t=(c_1,\dots,c_d)$. Then, clearly $\t\in I_{d,n}$, and
$\t<\t_i$, $i=1,2$. We have $w_i\le \t$, and hence $X_1\cap
X_2=X_{\t}$.

The proof when $X_1$ and $X_2$ are opposite Schubert varieties is
similar. The result for Richardson varieties follows immediately
from the result for Schubert varieties and the result for opposite
Schubert varieties.
\end{proof}

\begin{rem}\label{intersect}
Explicitly, in terms of the distributive lattice structure of
$I_{d, n}$, we have that $X_{w_1}^{v_1} \cap X_{w_2}^{v_2} =
X_{w_1 \wedge w_2}^{v_1 \vee v_2}$ (set theoretically), where
$w_1\wedge w_2$ is the {\it meet} of $w_1$ and $w_2$ (the largest
element of $W^{P_d}$ which is less than both $w_1$ and $w_2$) and
$v_1\vee v_2$ is the {\it join} of $v_1$ and $v_2$ (the smallest
element of $W^{P_d}$ which is greater than both $v_1$ and $v_2$).
The fact that $X_1 \cap X_2$ is reduced follows from \cite{m-r};
we will also provide a proof in Theorem \ref{reduced}.
\end{rem}

\section{Standard Monomial Theory for Richardson Varieties}
\label{smt1}

\subsection{Standard
Monomials}\label{sec:sm3.1} Let $R_0$ be the homogeneous
coordinate ring of $G_{d,n}$ for the Pl\"ucker embedding, and for
$w, v \in I_{d,n}$, let $R_w^v$ be the homogeneous coordinate ring
of the Richardson variety $X_w^v$. In this section, we present a
standard monomial theory for $X_w^v$  in the same spirit as in
\cite{Mu}. As mentioned in the introduction, standard monomial
theory consists in constructing an explicit basis for $R_w^v$.

\begin{defn}
A monomial $f=p_{\t_1}\cdots p_{\t_m}$ is said to be {\em
standard} if
\begin{equation*}
\t_1\ge \cdots \ge \t_m. \tag{*}
\end{equation*}
Such a monomial is said to be {\em standard on $X_w^v$}, if in
addition to condition (*), we have $w \ge \t_1$ and $\t_m \ge v$.
\end{defn}
\begin{rem}\label{nice}
Note that in the presence of condition (*), the standardness of
$f$ on $X_w^v$ is equivalent to the condition that $f\bigr|_
{X_w^v}\ne 0$. Thus given a standard monomial $f$, we have
$f\bigr|_ {X_w^v}$ is either $0$ or remains standard on $X_w^v$.
\end{rem}

\subsection{Linear Independence of Standard Monomials}
\begin{thm}\label{indep}
The standard monomials on $X_w^v$ of degree $m$ are linearly
independent in $R_w^v$.
\end{thm}
\begin{proof}
We proceed by induction on $\dim X_w^v$.

If $\dim X_w^v=0$, then $w = v$, $p_{w}^m$ is the only standard
monomial on $X_w^v$ of degree $m$, and the result is obvious. Let
$\dim X_w^v>0$. Let
\begin{equation*} 0 = \sum_{i=1}^rc_i F_i,\ c_i\in K^*\, ,\tag{$*$}
\end{equation*}
be a linear relation of standard monomials $F_i$ of degree $m$.
Let $F_i=p_{w_{i1}}\dots p_{w_{im}}$.  Suppose that $w_{i1}<w$ for
some $i$. For simplicity, assume that $w_{11}<w$, and $w_{11}$ is
a minimal element of $\{w_{j1}\mid w_{j1}<w\}$. Let us denote
$w_{11}$ by $\varphi$. Then for $i\ge 2$,
$F_i\res_{X_{\varphi}^v}$ is either $0$, or is standard on
$X_{\varphi}^v$. Hence restricting $(*)$ to $X_{\varphi}^v$, we
obtain a nontrivial standard sum on $X_{\varphi}^v$ being zero,
which is not possible (by induction hypothesis). Hence we conclude
that $w_{i1}=w$ for all $i$, $1\le i\le m$. Canceling $p_w$, we
obtain a linear relation among standard monomials on $X_w^v$ of
degree $m-1$. Using induction on $m$, the required result follows.
\end{proof}

\subsection{Generation by Standard Monomials}
\begin{thm}\label{generation}
Let $F=p_{w_1}\dots p_{w_m}$ be any monomial in the Pl\"ucker
coordinates of degree $m$. Then $F$ is a linear combination of
standard monomials of degree $m$.
\end{thm}
\begin{proof}
For $F=p_{w_1}\dots p_{w_m}$, define
$$N_F=l(w_1)N^{m-1}+l(w_2)N^{m-2}+\dots l(w_m),$$ where $ N\gg 0$,
say $N>d(n-d)$ ($=\dim G_{d,n}$) and $l(w)=\text{dim}X_w$. If $F$
is standard, there is nothing to prove. Let $t$ be the first
violation of standardness, i.e. $p_{w_1}\dots p_{w_{t-1}}$ is
standard, but $p_{w_1}\dots p_{w_t}$ is not. Hence $w_{t-1}\not\ge
w_t$, and using the quadratic relations (cf. Theorem \ref{T:greq})
\begin{equation*}
p_{w_{t-1}}p_{w_t}=\sum_{\a,\b}\pm p_\a p_\b,\tag{*}
\end{equation*}
$F$ can be expressed as $F=\sum F_i$, with $N_{F_i}>N_F$ (since
$\a>w_{t-1}$ for all $\a$ on the right hand side of $(*)$). Now
the required result is obtained by decreasing induction on $N_F$
(the starting point of induction, i.e. the case when $N_F$ is the
largest, corresponds to standard monomial $F=p_\te^m$, where
$\te=(n+1-d,n+2-d,\cdots ,n)$, in which case $F$ is clearly
standard).
\end{proof}
Combining Theorems \ref{indep} and \ref{generation}, we obtain
\begin{thm}\label{mainth}
Standard monomials on $X_w^v$ of degree $m$ give a basis for
$R_w^v$ of degree $m$.
\end{thm}

As a consequence of Theorem \ref{mainth} (or also Theorem
\ref{T:greq}), we have a qualitative description of a typical
quadratic relation on a Richardson variety $X_w^v$ as given by the
following
\begin{prop}\label{qualitative}
Let $w,\t,\varphi, v \in I_{d,n},\ w>\t,\varphi$ and $\t, \varphi
>v$. Further let $\t,\varphi$ be non-comparable (so that $p_\t
p_\varphi$ is a non-standard degree $2$ monomial on $X_w^v$). Let
\begin{equation*}
p_{\t}p_{\varphi}=\sum_{\a,\b}\ c_{\a,\b} p_\a p_\b,\ c_{\a,\b}\in
k^*\, \tag{*}
\end{equation*}
be the expression for $p_{\t}p_{\varphi}$ as a sum of standard
monomials on $X_w^v$. Then for every $\a,\b$ on the right hand
side we have, $\a>$ both $\t\text{ and }\varphi$, and $\b<$ both
$\t\text{ and }\varphi$.
\end{prop}
Such a relation as in (*) is called a {\em straightening
relation}.

\subsection{Equations Defining Richardson Varieties in the
Grassmannian}\label{defneqns} Let

\ni $w, v \in I_{d,n}$, with $w \ge v$. Let $\pi_w^v$ be the map
$R_0\to R_w^v$ (the restriction map). Let $\ker \pi_w^v=J_w^v$.
Let $Z_w^v=\{\text{all standard monomials }F\ |\ F \text{ contains
some } p_\varphi$ for some $w \not\ge \varphi \text{ or } \varphi
\not\ge v \}$. We shall now give a set of generators for $J_w^v$
in terms of Pl\"ucker coordinates.

\begin{lem}\label{ideal}
Let $I_w^v=(p_\varphi,w \not\ge \varphi \text{ or } \varphi
\not\ge v)$ (ideal in $R_0$). Then $Z_w^v$ is a basis for $I_w^v$.
\end{lem}

\begin{proof}
Let $F\in I_w^v$. Then writing $F$ as a linear combination of
standard monomials $$F=\sum a_i F_i +\sum b_j G_j,$$ where in the
first sum each $F_i$ contains some $p_{\t}$, with $w  \not\ge \t
\text{ or } \t \not\ge v $, and in the second sum each $G_j$
contains only coordinates of the form $p_{\t}$, with $w \ge \t \ge
v$. This implies that $\sum a_i F_i\in I_w^v$, and hence we obtain
$$\sum b_jG_j\in I_w^v.$$ This now implies that considered as an
element of $R_w^v$, $\sum b_jG_j$ is equal to $0$ (note that
$I_w^v \subset J_w^v$). Now the linear independence of standard
monomials on $X_w^v$ implies that $b_j=0$ for all $j$. The
required result now follows.
\end{proof}

\begin{prop}\label{eqns}
Let $w, v \in I_{d,n}$ with $w \ge v$. Then $R_w^v = R_0/I_w^v$.
\end{prop}
\begin{proof}
We have, $R_w^v = R_0/J_w^v$ (where $J_w^v$ is as above). We shall
now show that the inclusion $I_w^v \subset J_w^v$ is in fact an
equality. Let $F\in R_0$. Writing $F$ as a linear combination of
standard monomials $$F=\sum a_i F_i+\sum b_j G_j,$$ where in the
first sum each $F_i$ contains some term $p_{\t}$, with $w \not\ge
\t \text{ or } \t \not\ge v $, and in the second sum each $G_j$
contains only coordinates $p_{\t}$, with $w \ge \t \ge v$, we
have, $\sum a_i F_i\in I_w^v$, and hence we obtain

\ni $F\in J_w^v$

\ni $\Longleftrightarrow  \sum b_j G_j \in J_w^v$ (since $\sum a_i
F_i \in I_w^v$, and $I_w^v \subset J_w^v$)

\ni $ \Longleftrightarrow \pi_w^v(F)$ $(=\sum b_jG_j)$ is zero

\ni  $\Longleftrightarrow  \sum b_jG_j$ (= a sum of standard
monomials on $X_w^v$) is zero on $X_w^v$

\ni $\Longleftrightarrow b_j=0$ for all $j$ (in view of the linear
independence of standard monomials on $X_w^v$)

\ni $\Longleftrightarrow F=\sum a_i F_i$

\ni $\Longleftrightarrow F\in I_w^v$ .

\ni Hence we obtain $J_w^v=I_w^v$.
\end{proof}

\ni {\bf Equations defining Richardson varieties:}

Let $w, v \in I_{d,n}$, with $w \ge v$. By Lemma \ref{ideal} and
Proposition \ref{eqns}, we have that the kernel of $(R_0)_1\to
(R_w^v)_1$ has a basis given by $\{p_{\t} \mid w  \not\ge \t
\text{ or } \t \not\ge v \}$, and that the ideal $J_w^v$ (= the
kernel of the restriction map $R_0\to R_w^v$) is generated by
$\{p_{\t} \mid w \not\ge \t \text{ or } \t \not\ge v \}$. Hence
$J_w^v$ is generated by the kernel of $(R_0)_1\to (R_w^v)_1$. Thus
we obtain that $X_w^v$ is scheme-theoretically (even at the cone
level) the intersection of $G_{d,n}$ with all hyperplanes in
$\mathbb{P}(\wedge^dk^n)$ containing $X_w^v$. Further, as a closed
subvariety of $G_{d,n}$, $X_w^v$ is defined (scheme-theoretically)
by the vanishing of $\{ p_{\t} \mid w \not\ge \t \text{ or } \t
\not\ge v\}$.

\section[Standard Monomial Theory for a Union of Richardson
Varieties]{Standard Monomial Theory for a Union of Richardson
Varieties}\label{unn} In this section, we prove results similar to
Theorems \ref{indep} and \ref{mainth} for a union of Richardson
varieties.

Let $X_i$ be Richardson varieties in $G_{d,n}$. Let $X=\cup X_i$.

\begin{defn}
A monomial $F$ in the Pl\" ucker coordinates is {\em standard} on
the union $X=\cup X_i$ if it is standard on some $X_i$.
\end{defn}
\subsection{Linear Independence of Standard Monomials on $X=\cup
X_i$}

\begin{thm}\label{indep'}
Monomials standard on $X=\cup X_i$ are linearly independent.
\end{thm}

\begin{proof} If possible, let
\begin{equation*}
0 = \sum_{i=1}^r a_i F_i, \ a_i \in K^*\tag{*}
\end{equation*}
be a nontrivial relation among standard monomials on $X$. Suppose
$F_1$ is standard on $X_j$. Then restricting $(*)$ to $X_j$, we
obtain a nontrivial relation among standard monomials on $X_j$,
which is a contradiction (note that for any $i,\ F_i|_{X_j}$ is
either $0$ or remains standard on $X_j$; further, $F_1|_{X_j}$ is
non-zero).
\end{proof}

\subsection{Standard Monomial Basis}
\begin{thm}\label{union}
Let $X=\cup_{i=1}^r\,X_{w_i}^{v_i}$, and $S$ the homogeneous
coordinate ring of $X$. Then the standard monomials on $X$ give a
basis for $S$.
\end{thm}
\begin{proof}
For $w, v \in I_{d,n}$ with $w \ge v$, let $I_w^v$ be as in Lemma
\ref{ideal}. Let us denote $I_t=I_{w_t}^{v_t}, X_t=X_{w_t}^{v_t},
\, 1\le t\le r$. We have $R_{w_t}^{v_t}=R_0/I_t$ (cf. Proposition
\ref{eqns}). Let $S=R_0/I$. Then $I=\cap I_t$ (note that being the
intersection of radical ideals, $I$ is also a radical ideal, and
hence the set theoretic equality $X=\cup X_i$ is also scheme
theoretic).  A typical element in $R_0/I$ may be written as $\pi
(f)$, for some $f\in R_0$, where $\pi$ is the canonical projection
$R_0\to R_0/I$. Let us write $f$ as a sum of standard monomials
$$f=\sum a_j G_j+\sum b_lH_l,$$ where each $G_j$ contains some
$p_{\t_j}$ such that $w_i \not\ge \t_j \text { or } \t_j \not\ge
v_i$, for  $1\le i\le r$; and for each $H_l$, there is some $i_l$,
with $1\le i_l\le r$, such that $H_l$ is made up entirely of
$p_{\t}$'s with $w_{i_l} \ge \t \ge v_{i_l}$. We have $\pi(f)=\sum
b_l H_l$ (since $\sum a_j G_j\in I$). Thus we obtain that $S$ (as
a vector space) is generated by monomials standard on $X$. This
together with the linear independence of standard monomials on $X$
implies the required result.
\end{proof}

\subsection{Consequences}\label{consq}

\begin{thm}\label{reduced}
Let $X_1$, $X_2$ be two Richardson varieties in $G_{d,n}$. Then

$(1)$ $X_1\cup X_2$ is reduced.

$(2)$ If $X_1\cap X_2\ne\emptyset$, then $X_1\cap X_2$ is reduced.
\end{thm}
\begin{proof}
$(1)$ Assertion is obvious.

$(2)$ Let $X_1=X_{w_1}^{v_1}$, $X_2=X_{w_2}^{v_2}$,
$I_1=I_{w_1}^{v_1}$, and $I_2=I_{w_2}^{v_2}$. Let $A$ be the
homogeneous coordinate ring of $X_1\cap X_2$. Let $A=R_0/I$. Then
$I=I_1+I_2$. Let $F\in I$. Then by Lemma \ref{ideal} and
Proposition \ref{eqns}, in the expression for $F$ as a linear
combination of standard monomials $$F=\sum a_j F_j,$$ each $F_j$
contains some $p_\t$, where either ($(w_1 \text{ or } w_2) \not
\ge \t$) or ($\t \not\ge (v_1 \text{ or } v_2)$). Let $X_1\cap
X_2= X_{\m}^{\n}$ set theoretically, where $\m = w_1 \wedge w_2$
and $\n = v_1 \vee v_2$ (cf. Remark \ref{intersect}). If
$B=R_0/\sqrt{I}$, then by Lemma \ref{ideal} and Proposition
\ref{eqns}, under $\pi:R_0\to B$, $\ker\pi$ consists of all $f$
such that $f=\sum c_k f_k$, $f_k$ being standard monomials such
that each $f_k$ contains some $p_{\varphi}$, where $\m
\not\ge\varphi\text{ or } \varphi\not\ge \n$.  Hence either ($(w_1
\text{ or } w_2) \not \ge \varphi$) or ($\varphi \not\ge (v_1
\text{ or } v_2)$). Hence $\sqrt{I}=I$, and the required result
follows from this.
\end{proof}

\begin{defn}
Let $w > v$.  Define $\partial^+X_w^v := \bigcup\limits_{w > w'
\geq v}X_{w'}^v$, and $\partial^-X_w^v := \bigcup\limits_{w \geq
v' > v}X_w^{v'}$.
\end{defn}

\begin{thm}\label{pieri}(Pieri's formulas)  Let $w > v$.
\begin{enumerate}
\item $X_w^v\cap\{p_w=0\}= \partial^+ X_w^v,\quad\text{scheme
theoretically}.$

\item $X_w^v\cap\{p_v=0\}= \partial^- X_w^v,\quad\text{scheme
theoretically}.$

\end{enumerate}
\end{thm}
\begin{proof}
Let $X=\partial^+X_w^v$, and let $A$ be the homogeneous coordinate
ring of $X$. Let $A=R_w^v/I$. Clearly, $(p_w) \subseteq I,\ (p_w)$
being the principal ideal in $R_w^v$ generated by $p_w$. Let $f\in
I$. Writing $f$ as $$f=\sum b_i G_i+\sum c_j H_j,$$ where each
$G_i$ is a standard monomial in $R_w^v$ starting with $p_w$ and
each $H_j$ is a standard monomial in $R_w^v$ starting with
$p_{\te_{j1}}$, where $\te_{j1}<w$,  we have, $\sum b_i G_i\in I$.
This now implies $\sum c_j H_j$ is zero on $\partial^+X_w^v$. But
now $\sum c_j H_j$ being a sum of standard monomials on
$\partial^+X_w^v$, we have by Theorem \ref{indep'}, $c_j=0$, for
all $j$. Thus we obtain $f=\sum b_i G_i$, and hence $f\in (p_w)$.
This implies $I=(p_w)$. Hence we obtain $A=R_w^v/(p_w)$, and (1)
follows from this.  The proof of (2) is similar.
\end{proof}

\section{Vanishing Theorems}\label{van}

Let $X$ be a union of Richardson varieties. Let $S(X,m)$ be the
set of standard monomials on $X$ of degree $m$, and $s(X,m)$ the
cardinality of $S(X,m)$. If $X$=$X_w^v$ for some $w, v$, then
$S(X,m)$ and $s(X,m)$ will also be denoted by just $S(w, v ,m)$,
respectively $s(w, v, m)$.
\begin{lem}\label{count}
$(1)$ Let $Y=Y_1\cup Y_2$, where $Y_1$ and $Y_2$ are unions of
Richardson varieties such that $Y_1 \cap Y_2 \neq \emptyset$. Then
$$s(Y,m)=s(Y_1,m)+s(Y_2,m)-s(Y_1\cap Y_2,m).$$ $(2)$ Let $w > v$.
Then $$ \begin{aligned}
s(w,v,m)&=s(w,v,m-1)+s\left(\partial^+X_w^v,m\right)
\\ &= s(w,v,m-1)+s\left(\partial^-X_w^v,m\right)
\end{aligned}
$$.
\end{lem}
\ni (1) and (2) are easy consequences of the results of the
previous section.

Let $X$ be a closed subvariety of $G_{d,n}$.  Let
$L=p^*(\mathcal{O}_{\mathbb{P}}(1))$, where
$\mathbb{P}=\mathbb{P}(\wedge^dK^n)$, and
$p:X\hookrightarrow\mathbb{P}$ is the Pl\"ucker embedding
restricted to $X$.

\begin{prop}\label{inductive}
Let $r$ be an integer $\le d(n-d)$. Suppose that all Richardson
varieties $X$ in  $G_{d,n}$ of dimension at most $r$ satisfy the
following two conditions:

$(1)$ $H^i(X,L^m)=0$, for $i\ge 1$, $m\ge 0$.

$(2)$ The set $S(X,m)$ is a basis for $H^0(X,L^m),\ m\ge 0$.

\ni Then any union of Richardson varieties of dimension at most
$r$ which have nonempty intersection, and any nonempty
intersection of Richardson varieties, satisfy $(1)$ and $(2)$.
\end{prop}
\begin{proof}
The proof for intersections of Richardson varieties is clear,
since any nonempty intersection of Richardson varieties is itself
a Richardson variety (cf. Lemma \ref{inter} and Theorem
\ref{reduced}).

We will prove the result for unions by  induction on $r$. Let
$S_r$ denote the set of Richardson varieties $X$ in $G_{d,n}$ of
dimension at most $r$. Let $Y=\cup_{j=1}^t X_j$, $X_j\in S_r$. Let
$Y_1=\cup_{j=1}^{t-1} X_j$, and $Y_2=X_t$. Consider the exact
sequence $$0\to\mathcal{O}_Y\to\mathcal{O}_{Y_1}\oplus
\mathcal{O}_{Y_2}\to\mathcal{O}_{Y_1\cap Y_2}\to 0,$$ where
$\mathcal{O}_Y\to\mathcal{O}_{Y_1}\oplus \mathcal{O}_{Y_2}$ is the
map $f\mapsto (f\res_{Y_1},f\res_{Y_2})$ and
$\mathcal{O}_{Y_1}\oplus \mathcal{O}_{Y_2}\to\mathcal{O}_{Y_1\cap
Y_2}$ is the map $(f,g)\mapsto (f-g)\res_{Y_1\cap Y_2}$. Tensoring
with $L^m$, we obtain the long exact sequence $$\to
H^{i-1}(Y_1\cap Y_2,L^m)\to H^i(Y,L^m)\to H^i(Y_1,L^m)\oplus
H^i(Y_2,L^m)\to H^i(Y_1\cap Y_2,L^m)\to$$ Now $Y_1\cap Y_2$ is
reduced (cf. Theorem \ref{reduced}) and $Y_1\cap Y_2\in S_{r-1}$.
Hence, by the induction hypothesis $(1)$ and $(2)$ hold for
$Y_1\cap Y_2$. In particular, if $m\ge 0$, then $(2)$ implies that
the map $H^0(Y_1,L^m)\oplus H^0(Y_2,L^m)\to H^0(Y_1\cap Y_2,L^m)$
is surjective. Hence we obtain that the sequence $$0\to
H^0(Y,L^m)\to H^0(Y_1,L^m)\oplus H^0(Y_2,L^m)\to H^0(Y_1\cap Y_2,
L^m)\to 0$$ is exact. This implies $H^0(Y_1\cap Y_2, L^m)\to
H^1(Y,L^m)$ is the zero map; we have, $H^1(Y,L^m)\to
H^1(Y_1,L^m)\oplus H^1(Y_2,L^m)$ is also the zero map (since by
induction $H^1(Y_1,L^m)=0=H^1(Y_2,L^m)$). Hence we obtain
$H^1(Y,L^m)=0, \ m\ge 0$, and for $i\ge 2$, the assertion that
$H^i(Y,L^m)=0,\ m\ge 0$ follows from the long exact cohomology
sequence above (and induction hypothesis). This proves the
assertion (1) for $Y$.

To prove assertion (2) for $Y$, we observe $$
\begin{aligned}
h^0(Y,L^m)&=h^0(Y_1,L^m)+h^0(Y_2,L^m)-h^0(Y_1\cap Y_2,L^m) \\
&=s(Y_1,m)+s(Y_2,m)-s(Y_1\cap Y_2,m).
\end{aligned}
$$ Hence Lemma \ref{count} implies that $$h^0(Y,L^m)=s(Y,L^m).$$
This together with linear independence of standard monomials on
$Y$ proves assertion (2) for $Y$.
\end{proof}

\begin{thm}\label{main}
Let $X$ be a Richardson variety in $G_{d,n}$. Then

(a) $H^i(X,L^m)=0$ for $i\ge 1$, $m\ge 0$.

(b) $S(X,m)$ is a basis for $H^0(X,L^m)$, $m\ge 0$.
\end{thm}
\begin{proof}
We prove the result by induction on $m, \text{ and }\dim X$.

If $\dim X=0$, $X$ is just a point, and the result is obvious.
Assume now that $\dim X\ge 1$. Let $X=X_w^v, \, w > v$. Let
$Y=\partial^+{X_w^v}$. Then by Pieri's formula (cf. \S
\ref{pieri}), we have, $$Y=X(\t)\cap\{p_\t=0\}\quad\text{\rm
(scheme theoretically)}.$$ Hence  the  sequence
$$0\to\mathcal{O}_X(-1)\to\mathcal{O}_X\to\mathcal{O}_Y\to 0$$ is
exact. Tensoring it with $L^m$, and writing the cohomology exact
sequence, we obtain the long exact cohomology sequence $$\dots\to
H^{i-1}(Y,L^m)\to H^i(X,L^{m-1})\to H^i(X,L^m)\to
H^i(Y,L^m)\to\cdots.$$ Let $m\ge 0$, $i\ge 2$. Then the induction
hypothesis on $\dim X$ implies (in view of Proposition
\ref{inductive}) that $H^i(Y,L^{m})=0$, $i\ge 1$. Hence we obtain
that the sequence $0\to H^i(X,L^{m-1})\to H^i(X,L^{m})$, $i\ge 2$,
is exact. If $i=1$, again the induction hypothesis implies the
surjectivity of $H^0(X, L^m)\to H^0(Y, L^m)$. This in turn implies
that the map $H^0(Y,L^m)\to H^1(X,L^{m-1})$ is the zero map, and
hence we obtain that the sequence $0\to H^1(X, L^{m-1})\to
H^1(X,L^m)$ is exact. Thus we obtain that $0\to H^i(X,L^{m-1})\to
H^i(X,L^{m})$, $m\ge 0$, $i\ge 1$ is exact. But $H^i(X,L^m)=0$,
$m\gg 0$, $i\ge 1$ (cf. \cite{fac}). Hence we obtain
\begin{equation*}
H^i(X,L^m)=0 \text{ for } i\ge 1,\ m\ge 0,\tag{1}
\end{equation*}
and
\begin{equation*}
h^0(X,L^m)=h^0(X,L^{m-1})+h^0(Y,L^m).\tag{2}
\end{equation*}
In particular, assertion (a) follows from (1). The induction
hypothesis on $m$ implies that $h^0(X,L^{m-1})=s(X,m-1)$. On the
other hand, the induction hypothesis on $\dim X$ implies (in view
of Proposition \ref{inductive}) that $h^0(Y,L^m)=s(Y,m)$. Hence we
obtain
\begin{equation*}
h^0(X,L^m)=s(X,m-1)+s(Y,m).\tag{3}
\end{equation*}
Now $(3)$ together with Lemma \ref{count},(2) implies
$h^0(X,L^m)=s(X,m)$. Hence (b) follows in view of the linear
independence of standard monomials on $X_w^v$ (cf. Theorem
\ref{indep}).
\end{proof}

\begin{cor}\label{6.6.0.7}
We have

\ni 1. $R_w^v=\bigoplus\limits_{m\in \mathbb{Z}^+}\
H^0(X_w^v,L^{m}),\ w \geq v.$

\ni 2. $\dim H^0 (\partial^+ X_w^v, L^m) = \dim H^0 (\partial^-
X_w^v, L^m), \,w>v, \, m \ge 0.$
\end{cor}
\begin{proof}
Assertion 1 follows immediately from Theorems \ref{mainth} and
\ref{main}(b).  Assertion 2 follows from Lemma \ref{count},
Theorem \ref{inductive}(2), and Theorem \ref{main}(b).
\end{proof}

\section{Tangent Space and Smoothness}\label{tss}

\subsection{The Zariski Tangent Space} Let $x$ be a point on a
variety $X$. Let $\mathfrak{m}_x$ be the maximal ideal of the
local ring $\mathcal{O}_{X,\,x}$ with residue field
$K(x)$($=\mathcal{O}_{X,\,x}/\mathfrak{m}_x$). Note that $K(x)=K$
(since $K$ is algebraically closed). Recall that the Zariski
tangent space to $X$ at $x$ is defined as $$
\begin{aligned} T_x(X)&=\text{\rm
Der}_K(\mathcal{O}_{X,\,x},K(x))\\ &=\{D:\mathcal{O}_{X,\,x}\to
K(x), \ K\text{\rm -linear such that }D(ab)=D(a)b+aD(b)\}
\end{aligned} $$ (here $K(x)$ is regarded as an
$\mathcal{O}_{X,\,x}$-module). It can be seen easily that $T_x(X)$
is canonically isomorphic to $\text{\rm Hom}_{K\text{\rm
-mod}}(\mathfrak{m}_x/\mathfrak{m}_x^2,K)$.

\subsection{Smooth and Non-smooth Points}
A point $x$ on a variety $X$ is said to be a {\em simple} or {\em
smooth} or {\em nonsingular  point of $X$} if
$\mathcal{O}_{X,\,x}$ is a regular local ring. A point $x$ which
is not simple is called a {\em multiple} or {\em non-smooth} or
{\em singular point} of $X$. The set $\Sing X=\{x\in X\mid x\
\text{\rm  is a singular point}\}$ is called the {\em singular
locus of $X$}. A variety $X$ is said to be {\em smooth} if
$\text{\rm Sing}X=\emptyset$.  We recall the well known

\begin{thm}\label{1.2.37.1} \label{greater}
Let $x\in X$. Then $\dim_K T_x(X)\ge\dim\mathcal{O}_{X,\,x}$ with
equality if and only if $x$ is a simple point of $X$.
\end{thm}

\subsection{The Space $T_{w,\t}^v$} Let $G,T,B,P_d,W,R,W_{P_d}, R_{P_d}$
etc., be as in \S \ref{weylroot}.  We shall henceforth denote
$P_d$ by just $P$. For $\a \in R$, let $X_\a$ be the element of
the Chevalley basis for $\mathfrak{g} \ (={\rm Lie}G$),
corresponding to $\a$.  We follow \cite{bou} for denoting elements
of $R,R^{+}$ etc.

For $w \geq \t \geq v$, let $T_{w, \t}^v $ be the  Zariski tangent
space to $X_w^v$ at $e_{\t} $. Let $w_0$ be the element of largest
length in $W$. Now the tangent space to $G$ at $e_{id}$ is
$\mathfrak g$, and hence the tangent space to $G/P$ at $e_{id}$ is
$\bigoplus\limits_{\b\in R^+\setminus R_P^+}\  {{{\mathfrak
g}}}_{-\b}$. For $\t\in W$, identifying $G/P$ with $G/^{\t}P$
(where $^{\t}P=\t P\t ^{-1}$) via the map $gP\mapsto (n_\t
gn_\t^{-1}) {}^{\t}P,\ n_\t$ being a fixed lift of $\t$ in
$N_G(T)$, we have, the tangent space to $G/P$ at $e_\t$ is
$\bigoplus\limits_{\b
\in
\t(R^+)\setminus \t(R_P^+)}\  {{{\mathfrak g}}}_{-\b}$, i.e.,
$$T_{w_0,\t}^{\text{id}} =\bigoplus\limits_{\b
\in
\t(R^+)\setminus \t(R_P^+)}\  {{{\mathfrak g}}}_{-\b}. $$ Set
$$N_{w,\t}^v=\{ \b \in \t (R^+)\setminus \t(R_P^+)\ |\ X_{-\b }\in
T_{w,\t}^v \} .$$ Since $T_{w, \t}^v$ is a $T$-stable subspace of
$T_{w_0,\t}^v$, we have $$T_{w,\t}^v =\text{ the span of
}\{X_{-\b},\ \b \in N_{w,\t}^v\}.$$

\subsection{Certain Canonical Vectors in $T_{w,\t}^v$}
\label{curves}

For a root $\a \in R^+\setminus R_P^+$, let $Z_\a$ denote the
$SL(2)$-copy in $G$ corresponding to $\a$; note that $Z_\a$ is
simply the subgroup of $G$ generated by $U_\a$ and $U_{-\a}$.
Given $x\in W^P$, precisely one of $\{U_\a, U_{-\a}\}$ fixes the
point $e_x$. Thus $Z_\a\cdot e_x$ is a $T$-stable curve in $G/P$
(note that $Z_\a\cdot e_x\cong {\mathbb{P}}^1$), and conversely
any $T$-stable curve in $G/P$ is of this form (cf. \cite{C2}). Now
a $T$-stable curve $Z_\a\cdot e_x$ is contained in a Richardson
variety $X_w^v$ if and only if $e_x, e_{s_\a x}$ are both in
$X_w^v$.

\begin{lem}\label{curve}
Let $w,\t,v \in W^P,\ w\ge \t \ge v$. Let $\b \in \t(R^+\setminus
R_P^+)$. If $w\ge s_\b\t \ge v \, (\MMod W_P)$, then $X_{-\b} \in
T_{w,\t}^v$.
\end{lem}
\ni (Note that $s_{\b}\t$ need not be in $W^P$.)
\begin{proof}
The hypothesis that $w\ge s_\b\t \ge v (\MMod W_P)$ implies that
the curve $Z_\b\cdot e_\t$ is contained in $X_w^v$. Now the
tangent space to $Z_\b\cdot e_\t$ at $e_\t$ is the one-dimensional
span of $X_{-\b}$. The required result now follows.
\end{proof}
We shall show in Theorem \ref{tgt} that $w,\t,v$ being as above,
$T_{w,\t}^v$ is precisely the span of $\{X_{-\b},\ \b \in
\t(R^+\setminus R_P^+ )\ |\ w\ge s_\b\t \ge v \, (\MMod W_P) \}$.

\subsection{A Canonical Affine Neighborhood of a $T$-fixed Point
}\label{neigh}

Let $\t \in W$. Let $U^{-}_\t$ be the unipotent subgroup of $G$
generated by the root subgroups $U_{-\b},\ \b \in \t(R^+)$ (note
that $U^-_{\t}$ is the unipotent part of the Borel sub group
$^{\t}B^-$, opposite to $^{\t}B\ (= \t B\t^{-1})$). We have
$$U_{-\b} \cong \mathbb{G}_a,\ U^{-}_\t \cong \prod\limits_{\b \in
\t(R^+)}\ U_{-\b}. $$  Now, $U^-_{\t}$ acts on $G/P$ by left
multiplication. The isotropy subgroup in $U^-_{\t}$ at $e_\t$ is
$\Pi_{\b \in \t(R_P^+)}\ U_{-\b}$. Thus $U_{\t}^- e_\t \cong
\Pi_{\b \in \t(R^+\setminus R_P^+)}\ U_{-\b}$.  In this way,
$U^-_{\t} e_\t$ gets identified with $\mathbb{A}^{N}$, where $N=\#
(R^+\setminus R_P^+)$. We shall denote the induced coordinate
system on $U^{-}_\t e_\t$ by $\{x_{-\b},\ \b \in \t(R^+\setminus
R_P^+)\}$. In the sequel, we shall denote $U^- _{\t} e_{\t}$ by
${\mathcal O}^{-}_{\t}$ also. Thus we obtain that ${\mathcal
O}^{-}_{\t}$ is an affine neighborhood of $e_{\t}$ in $G/P$.

\subsection{The Affine Variety $Y_{w,\tau}^v$}\label{ywt} For $w,
\t, v \in W,\ w\ge \t \ge v$, let us denote $Y_{w,\t}^v:={\mathcal
O}^{-}_{\t}\cap X_{w, \t}^v$. It is a nonempty affine open
subvariety of $X_w^v$, and a closed subvariety of the affine space
${\mathcal O}^{-}_{\t}$.

Note that $L$, the ample generator of $\textrm{Pic}(G/P)$, is the
line bundle corresponding to the Pl\"ucker embedding, and
$H^0(G/P,L)=(\wedge^dK^n)^*$, which has a basis given by the
Pl\"ucker coordinates $ \{ p_{\te}, \, \te \in I_{d, n} \}$.  Note
also that the affine ring ${\mathcal O}^{-}_{\t}$ may be
identified as the homogeneous localization $(R_0)_{(p_{\t})}$,
$R_0$ being as in \S \ref{sec:sm3.1}. We shall denote
$p_{\te}/p_{\t}$ by $f_{\te,\t}$. Let $I_{w, \t}^v$ be the ideal
defining $Y_{w, \t}^v$ as a closed subvariety of ${\mathcal
O}^{-}_{\t}$. Then $I_{w, \t}^v$ is generated by $\{f_{\te, \t}
\mid w \not\geq \te \text{ or } \te \not\le v \}$.

\subsection{Basis for Tangent Space \& Criterion for Smoothness of $X_w^v$
at $e_\t$}

Let $Y$ be an affine variety in $\mathbb{A}^n$, and let $I(Y)$ be
the ideal defining $Y$ in $\mathbb{A}^{n}$.  Let $I(Y)$ be
generated by $\{f_1, f_2, \ldots, f_r\}$.  Let $J$ be the Jacobian
matrix $(\frac{\partial f_i }{\partial x_j})$.  We have (cf.
Theorem \ref{greater}) the dimension of the tangent space to $Y$
at a point $P$ is greater than or equal to the dimension of $Y$,
with equality if and only if $P$ is a smooth point; equivalently,
rank $J_P \leq \text{codim}_{\mathbb{A}^n}Y$ with equality if and
only if $P$ is a smooth point of $Y$ (here $J_P$ denotes $J$
evaluated at $P$).

Let $w,\t, v \in W,\ w\ge\t\ge v$. The problem of determining
whether or not  $e_{\tau}$ is a smooth point of $X_w^v$ is
equivalent to determining whether or not $e_{\tau}$ is a smooth
point of $Y_{w,\tau}^v$ (since $Y_{w,\t}^v$ is an open
neighborhood of $e_{\tau}$ in $X_w^v$). In view of Jacobian
criterion, the problem is reduced to computing $(\partial f_{\te,
\t} / \partial x_{-\b})_{e_\t}, \, w \not\ge \te \text{ or } \te
\not\ge v$ (the Jacobian matrix evaluated at $e_\tau$). To carry
out this computation, we first observe the following:

Let $V$ be the $G$-module $H^0(G/P,L) \ ( = (\wedge^dK^n)^*)$. Now
$V$ is also a ${\mathfrak g}$-module. Given $X$ in ${\mathfrak
g}$, we identify $X$ with the corresponding right invariant vector
field $D_X$ on $G$. Thus we have $D_X p_{\te} = X p_{\te}$, and we
note that
\begin{equation*}
(\partial f_{\te, \t} / \partial x_{-\b})({e_\t}) =
X_{-\b}\,p_{\te} (e_\t), \ \b \in \t(R^+\setminus R_P^+),
\end{equation*}
where the left hand side denotes the partial derivative evaluated
at $e_\t$.

We make the folowing three observations:

(1) For $\te,\mu\in W^P$, $p_\te(e_\mu)\ne 0\iff\te=\mu$, where,
recall that for $\te =(i_1\cdots i_d) \in W^P,\,\ e_\te$ denotes
the vector $e_{i_1}\wedge \cdots \wedge e_{i_d}$ in $\wedge^d
K^n$, and $p_\te$ denotes the Pl\"ucker coordinate associated to
$\te$.

(2) Let $X_\a$ be the element of the Chevalley basis of
$\mathfrak{g}$, corresponding to $\a \in R$. If $X_\a p_\mu\ne 0$,
$\mu\in W^P$, then $X_\a p_\mu=\pm p_{s_\a \mu}$, where $s_\a$ is
the reflection corresponding to the root $\a$.

(3) For $\a\ne \b$, if $X_\a p_\mu,\ X_\b p_\mu$ are non-zero,
then $X_\a p_\mu \ne X_\b p_\mu$.

The first remark is obvious, since $\{p_\te\mid\te\in W^P\}$ is
the basis of $(\wedge^P K^n)^* \ (=H^0(G/P,L))$, dual to the basis
$\{e_\varphi,\ \varphi \in I_{d,n}\}$ of $\wedge^d K^n$. The
second remark is a consequence of $SL_2$ theory, using the
following facts:

(a)
$|\langle\chi,\a^*\rangle|=\left|\frac{2(\chi,\a)}{(\a,\a)}\right|=0\
\text{\rm or } 1$, $\chi$ being the weight of $p_\mu$.

(b) $p_\mu$ is the lowest weight vector for the Borel subgroup
${}^\mu B=\mu B\mu^{-1}$.

The third remark follows from weight considerations (note that if
$X_\a p_\mu \ne 0$, then $X_\a p_\mu$ is a weight vector (for the
$T$-action) of weight $\chi +\a$, $\chi$ being the weight of
$p_\mu$).

\begin{thm}\label{jacobian}
Let $w, \t, v \in W^P,\ w \ge \t \ge v$. Then $$\dim
T_{w,\t}^v=\#\{\g\in \t(R^+\setminus R_P^+)\mid w\ge s_\g \t \ge v
\, (\MMod W_P)\}.$$
\end{thm}
\begin{proof} By Lemma \ref{ideal} and Proposition \ref{eqns}, we have,
$I_{w,\t}^v$ is generated by $\{f_{\te, \t} \mid w \not\ge\te
\text{ or } \te \not \ge v \}$. Denoting the affine coordinates on
$\mathcal{O}^{-}_\t $ by $x_{-\b},\ \b \in \t(R^+\setminus
R_P^+)$, we have the evaluations of $\frac{\partial f_{\te,
\t}}{\partial x_\b}$ and $X_\b p_\te$ at $e_\t$ coincide. Let
$J_w^v$ denote the Jacobian matrix of $Y_{w,\t}^v$ (considered as
a subvariety of the affine space $\mathcal{O}^{-}_\t $). We shall
index the rows of $J_w^v$ by $\{f_{\te, \t} \mid w \not\ge\te
\text{ or } \te \not \ge v \}$ and the columns by $x_{-\b},\ \b
\in \t(R^+\setminus R_P^+)$. Let $J_w^v(\t)$ denote $J_w^v$
evaluated at $e_\t$. Now in view of (1) \& (2) above, the
$(f_{\te,\t},x_{-\b})$-th entry in $J_w^v(\t)$ is non-zero if and
only if $X_\b p_\te=\pm p_\t$. Hence in view of (3) above, we
obtain that in each row of $J_w^v(\t)$, there is at most one
non-zero entry. Hence rank$J_w^v(\t)=$ the number of non-zero
columns of $J_w^v(\t)$. Now, $X_\b p_\te=\pm p_\t$ if and only if
$\te \equiv s_\b \t \, (\MMod W_P)$. Thus the column of
$J_w^v(\t)$ indexed by $x_{-\b}$ is non-zero if and only if
$w\not\ge s_\b \t \, (\MMod W_P)$ or $s_\b \t \not\ge v \, (\MMod
W_P) $. Hence rank$J_w^v(\t)=\#\{\g \in \t(R^+\setminus R_P^+)\ |\
w\not\ge s_\g \t \, (\MMod W_P) \text{ or } s_\g \t \not\ge v \,
(\MMod W_P)\}$ and thus we obtain $$\dim T_{w,\t}^v=\#\{\g\in
\t(R^+\setminus R_P^+)\mid w\ge s_\g \t \ge v \, (\MMod W_P) \}.$$
\end{proof}

\begin{thm}\label{tgt}
Let $w, \t, v$ be as in Theorem \ref{jacobian}. Then $\{X_{-\b},\
\b \in \t(R^+\setminus R_P^+)\ |\ w\ge s_\b \t \ge v \, (\MMod
W_P)\}$ is a basis for $T_{w,\t}^v$.
\end{thm}
\begin{proof}
Let $\b \in \t(R^+\setminus R_P^+)$ be such that $w\ge s_\b \t \ge
v \, (\MMod W_P)$. We have (by Lemma \ref{curve}), $X_{-\b} \in
T_{w,\t}^v$. On the other hand, by Theorem \ref{jacobian},

\ni $\dim T_{w,\t}^v=\#\{\b\in \t(R^+\setminus R_P^+)\mid w\ge
s_\b \t \ge v \, (\MMod W_P)\}.$ The result follows from this.
\end{proof}

\begin{cor}\label{smcrit1}
$X_w^v$ is smooth at $e_\t$ if and only if $l(w)-l(v)=\#\{\a \in
R^+\setminus R_P^+ \mid w\ge \t s_\a \ge v \, (\MMod W_P) \}$.
\end{cor}
\begin{proof}
We have, $X_w^v$ is smooth at $e_\t$ if and only if $\dim
T_{w,\t}^v=\dim X_w^v$, and the result follows in view of
Corollary \ref{morerich} and Theorem \ref{jacobian} (note that if
$\b=\t(\a)$, then $s_\b \t=\t s_\a \, (\MMod W_P)$).
\end{proof}

\section{Multiplicity at a Singular Point}\label{mas}

\subsection{Multiplicity of an Algebraic Variety at a Point}

Let $B$ be a graded, affine $K$-algebra such that $B_1$ generates
$B$ (as a $K$-algebra). Let $X = \textrm{Proj}(B)$. The function
$h_B(m)$ (or $h_X(m)$) $= \textrm{dim}_KB_m$, $m \in \mathbb{Z}$
is called the {\it Hilbert function} of $B \,(\textrm{or } X)$.
There exists a polynomial $P_B(x)$ (or $P_X(x)$) $\in
\mathbb{Q}[x]$, called the {\it Hilbert polynomial} of $B\,
(\textrm{or } X)$, such that $f_B(m) = P_B(m)$ for $m \gg 0$. Let
$r$ denote the degree of $P_B(x)$. Then $r = \textrm{dim}(X)$, and
the leading coefficient of $P_B(x)$ is of the form $c_B/r!$, where
$c_B \in \mathbb{N}$. The integer $c_B$ is called the {\it degree
of $X$}, and denoted $\textrm{deg}(X)$ (see \cite{eis} for
details). In the sequel we shall also denote $\textrm{deg}(X)$ by
$\textrm{deg}(B)$.

Let $X$ be an algebraic variety, and let $P \in X$.  Let $A =
\mathcal{O}_{X, P}$ be the stalk at $P$ and $\mathfrak{m}$ the
unique maximal ideal of the local ring $A$. Then the {\it tangent
cone} to $X$ at $P$, denoted $\textrm{TC}_P(X)$, is
$\textrm{Spec}(\textrm{gr}(A,\mathfrak{m}))$, where
$\textrm{gr}(A,\mathfrak{m}) =
\oplus_{j=0}^{\infty}\mathfrak{m}^j/\mathfrak{m}^{j+1}$. The {\it
multiplicity} of $X$ at $P$, denoted $\Mult_P(X)$, is
$\textrm{deg}(\textrm{Proj}(\textrm{gr}(A,\mathfrak{m})))$. (If $X
\subset K^n$ is an affine closed subvariety, and $m_P \subset
K[X]$ is the maximal ideal corresponding to $P \in X$, then
$\textrm{gr}(K[X], m_P) = \textrm{gr}(A,\mathfrak{m})$.)

\subsection{Evaluation of Pl\" ucker Coordinates on
$U_\t^-e_\t$}\label{evpl} Let $X=X_w^v$. Consider a $\t\in W^P$
such that $w \ge \t \ge v$.

\ni I. Let us first consider the case $\t=\text{\rm id}$. We
identify $U^-e_{\text{\rm id}}$ with $$ \left\{
\begin{pmatrix} &\text{\rm Id}_{d\times d}\\ x_{d+1\, 1} &\dots &x_{d+1\,
d}\\ \vdots &&\vdots\\ x_{n\, 1}&\dots&x_{n\,d}
\end{pmatrix},\quad x_{ij}\in k,\quad d+1\le i\le n, 1\le j\le d
\right\}. $$

Let $A$ be the affine algebra of $U^-e_{\text{\rm id}}$. Let us
identify $A$ with the polynomial algebra $k[x_{-\b},\b\in
R^+\setminus R_P^+]$. To be very precise, we have
$R^+\,\setminus\, R_P^+ =\{\e_j-\e_i$, $1\le j\le d,\  d+1\le i\le
n\}$; given $\b\in R^+\setminus R_P^+$, say $\b=\e_j-\e_i$, we
identify $x_{-\b}$ with $x_{ij}$. Hence we obtain that the
expression for $f_{\te,\textrm{id}}$ in the local coordinates
$x_{-\b}$'s is homogeneous.
\begin{ex} Consider $G_{2,4}$. Then
$$ U^-e_{\text{\rm id}}= \left\{
\begin{pmatrix} 1&0\\ 0&1\\ x_{31}&x_{32}\\ x_{41}&x_{42}
\end{pmatrix},\quad x_{ij}\in k
\right\}. $$ On $U^-e_{\text{\rm id}}$, we have $p_{12}=1$,
$p_{13}=x_{32}$, $p_{14}=x_{42}$, $p_{23}=x_{31}$,
$p_{24}=x_{41}$, $p_{34}=x_{31}x_{42}-x_{41}x_{32}$.
\end{ex}
Thus a Pl\"ucker coordinate is homogeneous in the local
coordinates $x_{ij},\ d+1\le i\le n,\ 1\le j\le d$.

\ni II. Let now $\t$ be any other element in $W^P$, say
$\t=(a_1,\dots,a_n)$. Then $U_\t^-e_\t$ consists of $\{N_{d,n}\}$,
where $N_{d,n}$ is obtained from $\begin{pmatrix}\text{\rm
Id}\\X\end{pmatrix}_{n\times d}$ (with notations as above) by
permuting the rows by $\t^{-1}$. (Note that $U_\t^-e_\t=\t
U^-e_{\text{\rm id}}$.)

\begin{ex} Consider $G_{2,4}$, and  let $\t=(2314)$. Then $\t^{-1}=(3124)$,
and $$U_\t^-e_\t= \left\{
\begin{pmatrix} x_{31}&x_{32}\\ 1&0\\ 0&1\\ x_{41}&x_{42}
\end{pmatrix},\quad x_{ij}\in k
\right\}. $$

We have on $U_\t^-e_\t$, $p_{12}=-x_{32}$, $p_{13}=x_{31}$,
$p_{14}=x_{31}x_{42}-x_{41}x_{32}$, $p_{23}=1$, $p_{24}=x_{42}$,
$p_{34}=-x_{41}$.
\end{ex}

As in the case $\t=\text{\rm id}$, we find that for $\te\in W^P$,
$f_{\te,\t}:=$ $p_\te\res_{U_\t^-e_\t}$ is homogeneous in local
coordinates. In fact we have
\begin{prop} Let $\te\in W^P$. We have a natural isomorphism
$$k[x_{-\b},\b\in R^+\setminus R_P^+]\cong k[x_{-\t(\b)}, \b\in
R^+\setminus R_P^+],$$ given by $$f_{\te,\textrm{id}}\mapsto
f_{\t\te,\t}.$$
\end{prop}

\ni The proof is immediate from the above identifications of
$U^-e_{\text{\rm id}}$ and $U_\t^-e_\t$. As a consequence, we have
\begin{cor}\label{8.8.3.5} Let $\te\in W^P$. Then the polynomial expression
for $f_{\te,\t}$ in the local coordinates $\{x_{-\t(\b)}, \b\in
R^+\setminus R_P^+\}$ is homogeneous.
\end{cor}

\subsection{The algebra $A_{w,\t}^v$}\label{8.8.4}

As above, we identify $A_\t$, the affine algebra of $U_\t^-e_\t$
with the polynomial algebra $K[x_{-\b},\b\in\t(R^+\setminus
R_P^+)]$. Let $A_{w,\t}^v=A_\t/I_{w,\t}^v$, where $I_{w\t}^v$ is
the ideal of elements of $A_\t$ that vanish on $X_w^v \cap
U_\t^-e_\t$.

Now $I(X_w^v)$, the ideal of $X_w^v$ in $G/P$, is generated by
$\{p_\te,\te\in W^P\mid w\not\ge\te \text{ or } \te \not\ge v\}$.
Hence we obtain (cf. Corollary \ref{8.8.3.5}) that $I_{w,\t}^v$ is
homogeneous. Hence we get
\begin{equation*}
\text{gr\,}(A_{w,\t}^v,M_{w,\t}^v)=A_{w,\t}^v,\tag{*}
\end{equation*}
where $M_{w,\t}^v$ is the maximal ideal of $A_{w,\t}^v$
corresponding to $e_\t$. In particular, denoting the image of
$x_{-\b}$ under the canonical map $A_\t\to A_{w,\t}^v$ by just
$x_{-\b}$, the set $\{x_{-\b}\mid\b\in \t(R^+\setminus R_P^+)\}$
generates $A_{w,\t}^v$. Let $R_w^v$ be the homogeneous coordinate
ring of $X_w^v$ (for the Pl\" ucker embedding),
$Y_{w,\t}^v=X_w^v\cap U_\t^-e_\t.$ Then $K[Y_{w,
\t}^v]=A_{w,\t}^v$ gets identified with the homogeneous
localization ${(R_w^v)}_{(p_\t)}$, i.e. the subring of
${(R_w^v)}_{p_\t}$ (the localization of $R_w^v$ with respect to
$p_\t$) generated by the elements $$\{\frac{p_\te}{p_\t},\te\in
W^P,w \ge \te \ge v\}.$$

\subsection{The Integer $\deg_{\t}(\te)$} Let $\te\in W^P$. We
define $\deg_{\t}(\te)$ by $$\deg_{\t}(\te):=\deg f_{\te,\t}$$
(note that $f_{\te,\t}$ is homogeneous, cf. Corollary
\ref{8.8.3.5}). In fact, we have an explicit expression for
$\deg_{\t}(\te)$, as follows (cf. \cite{L-G}):
\begin{prop} Let $\te\in W^P$. Let $\t=(a_1,\dots,a_n)$, \
$\te=(b_1,\dots,b_n)$. Let
$r=\#\{a_1,\dots,a_d\}\cap\{b_1,\dots,b_d\}$. Then
$\deg_{\t}(\te)=d-r$.
\end{prop}

\subsection{A Basis for the Tangent Cone}

Let $Z_\t=\{\te \in W^P \,|\,\text{either }\te\ge \t\,\text{or
}\t\ge \te\}.$

\begin{thm} With notations as above, given $r\in\mathbb{Z}^+$,

$$\{f_{\te_1,\t}\dots f_{\te_m,\t}\mid w\ge\te_1\ge...\ge\te_m\ge
v ,\, \te_i\in Z_\t,\ \sum_{i=1}^m\deg_{\t}(\te_i)=r \}$$ is a
basis for $(M_{w,\t}^v)^r/(M_{w,\t}^v)^{r+1}$.
\end{thm}
\begin{proof} For  $F=p_{\te_1}\cdots p_{\te_m}$, let
deg$F$ denote the degree of $f_{{\te_1},\t}\cdots f_{{\te_m},\t}$.
Let $$A_r=\{F=p_{\te_1}\cdots p_{\te_m},\ w\ge \te_i\ge v
\,|\,\text{deg}F=r\}.$$ Then in view of the relation (*) in \S
\ref{8.8.4}, we have, $A_r$ generates
$(M^v_{w,\t})^r/(M^{v}_{w,\t})^{r+1}$. Let $F\in A_r$, say
$F=p_{\t_1}\cdots p_{\t_m}$. From the results in \S \ref{smt1}, we
know that $p_{\t_1}\cdots p_{\t_m}$ is a linear combination of
standard monomials $p_{\te_1}\cdots p_{\te_m},\ w\ge \te_i \ge v$.
We {\em claim} that in each $p_{\te_1}\cdots p_{\te_m},\ \te_i \in
Z_\t, \text{ for all }i$. Suppose that for some $i,\te_i \not\in
Z_\t$. This means $\te_i$ and $\t$ are not comparable. Then using
the fact that $f_{\t,\t}=1$, on $Y_{w,\t}^v$, we replace
$p_{\te_i}$ by $p_{\te_i}p_\t$ in $p_{\te_1}\cdots p_{\te_m}$. We
now use the straightening relation (cf. Proposition
\ref{qualitative}) $p_{\te_i}p_\t=\sum c_{\a,\b}p_\a p_\b$ on
$X_w^v$, where in each term $p_\a p_\b$ on the right hand side, we
have $\a<w$, and $\a
>$ both $\te_i$ and $\t$, and $\b <$ both $\te_i$ and $\t$; in
particular, we have, in each term $p_\a p_\b$ on the right hand
side, $\a, \b$ belong to $Z_{\t}$. We now proceed as in the proof
of Theorem \ref{generation} to conclude that on $Y_{w,\t}^v,\ F$
is a linear combination of standard monomials in $p_\te$'s, $\te
\in Z_{\t}$ which proves the Claim.

Clearly, $\{f_{\te_1,\t}\dots f_{\te_m,\t}\mid
w\ge\te_1\ge...\ge\te_m\ge v  , \, \te_i\in Z_\t \}$ is linearly
independent in view of Theorem \ref{indep} (Since $p_{\t}^l
f_{\te_1,\t} \cdots f_{\te_m,\t} = p_{\t}^{l-m}p_{\te_1} \cdots
p_{\te_m}$ for $l\ge m$, and the monomial on the right hand side
is standard since $\te_i \in Z_{\t}$).
\end{proof}

\subsection{Recursive Formulas for $\Mult_\t{X_w^v}$}

\begin{defn}
If $w > \t \geq v$, define $\partial_{w,\t}^{v, +} := \{ w' \in
W^P \mid w > w' \geq \t \geq v, \, l(w')=l(w)-1\}$.  If $w \geq \t
>v$, define $\partial_{w,\t}^{v, -} := \{ v' \in W^P \mid w \geq
\t \geq v' > v, l(v')=l(v)+1\}$.
\end{defn}
\begin{thm}\label{8.8.5.3}
\begin{enumerate}
\item  Suppose $w > \t \ge v$. Then
$$(\Mult_\t{X_w^v})\deg_\t w=\sum\limits_ {w' \in \partial_{w,
\t}^{v, +}} \Mult_\t X_{w'}^v.$$

\item   Suppose $w \ge \t > v$. Then
$$(\Mult_\t{X_w^v})\deg_\t v=\sum\limits_ {v' \in \partial_{w,
\t}^{v, -}} \Mult_\t X_{w}^{v'}.$$

\item $\Mult_\t{X_\t^\t} = 1$.

\end{enumerate}
\end{thm}
\begin{proof}
Since $X_{\t}^{\t}$ is a single point, (3) is trivial.  We will
prove (1); the proof of (2) is similar.

Let $H_{\t} = \bigcup\limits_{w' \in \partial_{w, \t}^{v, +
}}X_{w'}^v$. Let $\varphi_w^v(r)$ (resp. $\varphi_{H_{\t}}(r)$) be
the Hilbert function for the tangent cone of $X_w^v$ (resp.
$H_{\t}$) at $e(\t)$, i.e. $$\varphi_w^v(r)=\dim
((M_{w,\t}^v)^r/(M_{w,\t}^v)^{r+1}).$$ Let $$
\mathcal{B}_{w,\t}^v(r)=\left\{p_{\t_1}\dots p_{\t_m},\,\t_i\in
Z_{\t}\mid (1)\, w\ge\t_1\ge\dots\ge\t_m\ge v, (2) \sum
\deg_{\t}(\t_i)=r\right\}. $$ Let $$
\begin{aligned}
\mathcal{B}_1&=\{p_{\t_1}\dots
p_{\t_m}\in\mathcal{B}_{w,\t}^v(r)\mid\t_1=w\},\\
\mathcal{B}_2&=\{p_{\t_1}\dots
p_{\t_m}\in\mathcal{B}_{w,\t}^v(r)\mid\t_1<w\}.
\end{aligned}
$$ We have $\mathcal{B}_{w, \t
}^v(r)=\mathcal{B}_1\dot\cup\mathcal{B}_2$. Hence denoting
$\deg_{\t}(w)$ by $d$, we obtain
$$\varphi_w^v(r+d)=\varphi_w^v(r)+\varphi_{H_\t}(r+d).$$ Taking
$r\gg 0$ and comparing the coefficients of $r^{u-1}$, where
$u=\dim X_w^v$, we obtain the result.
\end{proof}

\begin{cor} Let $w > \t > v$.  Then
$$(\Mult_{\t}{X_w^v}) (\deg_{\t}w + \deg_{\t}v) = \sum\limits_ {w'
\in \partial_{w, \t}^{v, +}} \Mult_\t{X_{w'}^v} + \sum\limits_ {v'
\in\partial_{w, \t}^{v,- } }\Mult_\t{X_w^{v'}}.$$
\end{cor}

\begin{thm}\label{prodmult}  Let $w \ge \t \ge v$.  Then
$\Mult_{\t}{X_w^v} = (\Mult_\t{X_w})\cdot(\Mult_\t{X^v})$
\end{thm}
\begin{proof}
We proceed by induction on $\dim X_w^v$.

If $\dim X_w^v = 0$, then $w = \t = v$.  In this case, by Theorem
{\ref{8.8.5.3}} (3), we have that $\Mult_\t{X_\t^\t} = 1$. Since
$e_{\t} \in Be_{\t } \subseteq X_w$, and $Be_{\t}$ is an affine
space open in $X_w$, $e_{\t}$ is a smooth point of $X_w$, i.e.
$\Mult_\t{X_w} = 1$.  Similarly, $\Mult_\t{X^v} = 1$.

Next suppose that $\dim X_w^v > 0$, and $w > \t \ge v$.  By
Theorem {\ref{8.8.5.3}} (1), $$\begin{aligned} \Mult_\t{X_w^v} &=
\frac{1}{\deg_{\t}w}\sum\limits_ {w' \in
{\partial_{w,\t}^{v,+}}}\Mult_\t{X_{w'}^v}
\\ &= \frac{1}{\deg_{\t}w}\sum\limits_
{w' \in {\partial_{w,\t}^{v,+}}}\Mult_\t{X_{w'}}\cdot
\Mult_\t{X^v}
\\ &= \left(\frac{1}{\deg_{\t}w}\sum\limits_
{w' \in {\partial_{w,\t}^{v,+}}}\Mult_\t{X_{w'}}\right)\cdot
\Mult_\t{X^v}
\\ &= \left(\frac{1}{\deg_{\t}w}\sum\limits_
{w' \in {\partial_{w,\t}^{v,+}}}\Mult_{\t}{X_{w'}^
{\text{id}}}\right)\cdot \Mult_\t{X^v}
\\ &= \left(\Mult_{\t}{X_w^{\text{id}}}\right)\cdot \Mult_\t{X^v}
= \Mult_\t{X_w} \cdot \Mult_\t{X^v}.
\end{aligned}$$

The case of $\dim X_w^v > 0$ and $w = \t > v$ is proven similarly.
\end{proof}

\begin{cor}\label{smcrit2}
Let $w \ge \t \ge v$.  Then $X_w^v$ is smooth at $e_\t$ if and
only if both $X_w$ and $X^v$ are smooth at $e_\t$.
\end{cor}

\begin{rem}
The following alternate proof of Theorem \ref{prodmult} is due to
the referee, and we thank the referee for the same.

\vskip.3cm
\begin{quote}{\it
Identify $\mathcal{O}^-_{\tau}$ with the affine space $\AAA ^N$
where $N=d(n-d)$, by the coordinate functions defined in \S
\ref{neigh}. Then $X_w\cap \mathcal{O}^-_{\tau}$ and $X^v\cap
\mathcal{O}^-_{\tau}$ are closed subvarieties of that affine
space, both invariant under scalar multiplication (e.g. by
Corollary \ref{8.8.3.5}). Moreover, $X_w$ and $X^v$ intersect
properly along the irreducible subvariety $X_w^v$; in addition,
the Schubert cells $C_w$ and $C^v$ intersect transversally (by
\cite{Ri}).

Now let $Y$ and $Z$ be subvarieties of $\AAA^N$, both invariant
under scalar multiplication, and intersecting properly. Assume in
addition that they intersect transversally along a dense open
subset of $Y\cap Z$. Then $$ \Mult_{\bf{o}}(Y\cap
Z)=\Mult_{\bf{o}}(Y) \cdot \Mult_{\bf{o}}(Z) $$ where ${\bf{o}}$
is the origin of ${\AAA ^N}$.

To see this, let $\PPP(Y)$, $\PPP(Z)$ be the closed subvarieties
of $\PPP(\AAA^N)=\PPP^{N-1}$ associated with $Y$, $Z$. Then
$\Mult_{\bf{o}}(Y)$ equals the degree ${\rm deg}(\PPP(Y))$, and
likewise for $Z$, $Y\cap Z$. Now $$ {\rm deg}(\PPP(Y\cap Z))= {\rm
deg}(\PPP (Y)\cap\PPP(Z)) ={\rm deg}(\PPP(Y))\cdot{\rm
deg}(\PPP(Z)) $$ by the assumptions and the Bezout theorem (see
\cite{ful}, Proposition 8.4 and Example 8.1.11).}
\end{quote}

\vskip.3cm

It has also been pointed out by the referee that the above
alternate proof in fact holds for Richardson varieties in a
minuscule $G/P$, since the intersections of Schubert and opposite
Schubert varieties with the opposite cell are again invariant
under scalar multiplication (the result analogous to Corollary
\ref{8.8.3.5} for a minuscule $G/P$ follows from the results in
\cite{l-w}). Recall that for $G$ a semisimple algebraic group and
$P$ a maximal parabolic subgroup of $G$, $G/P$ is said to be {\em
minuscule} if the associated fundamental weight $\o$ of $P$
satisfies $$(\o,\b^*)\, (\, =2(\o,\b)/(\o,\b))\, \leq 1$$ for all
positive roots $\b$, where $(\ , \, )$ denotes a $W$-invariant
inner product on $X(T)$.

\end{rem}

\subsection{Determinantal Formula for $\Mult X_w^v$}

In this section, we extend the

\ni Rosenthal-Zelevinsky determinantal formula (cf.
\cite{rose-zel}) for the multplicity of a Schubert variety at a
$T$-fixed point to the case of Richardson varieties.  We use the
convention that the binomial coefficient ${a \choose b} = 0$ if
$b<0$.

\begin{thm}\label{ros-zel} {(Rosenthal-Zelevinsky)} Let $w =
(i_1,\ldots,i_d)$ and $\t=(\t_1,\ldots,\t_d)$ be such that
$w\ge\t$. Then

$$\Mult_{\tau}{X_w} = (-1)^{\kappa_1 + \cdots + \kappa_d} \left|
\begin{matrix}{i_1 \choose -\kappa_1} & \cdots & {i_d \choose
-\kappa_d} \cr {i_1 \choose 1-\kappa_1} & \cdots & {i_d \choose
1-\kappa_d} \cr \vdots & & \vdots \cr {i_1 \choose d-1-\kappa_1} &
\cdots & {i_d \choose d-1-\kappa_d}
\end{matrix}
\right|,$$ where $\kappa_q := $ $\# \{ \t_p \mid \t_p > i_q \}$,
for $q = 1, \ldots, d$.
\end{thm}

\begin{lem}\label{oppreg}
$\Mult_{\t}{X^v} = \Mult_{w_o\t}{X_{w_0v}}$, where $w_0 =
(n+1-d,\ldots,n)$.
\end{lem}
\begin{proof}
Fix a lift $n_0$ in $N(T)$ of $w_0$.  The map $f:X^v\to n_0X^v$
given by left multiplication is an isomorphism of algebraic
varieties. We have $f(e_\t)=e_{w_0\t}$, and
$n_0X^v=n_0\overline{B^-e_v}$
$=n_0\overline{n_0Bn_0e_v}=\overline{Bn_0e_v}$
$=\overline{Be_{w_0v}}=X_{w_0v}$.
\end{proof}

\begin{thm}\label{rzr} Let $w = (i_1,\ldots,i_d)$,
$\t=(\t_1,\ldots,\t_d)$, and $v=(j_1,\ldots,j_d)$ be such that
$w\ge\t\ge v$.  Then

$$\Mult_{\tau}{X_w^v} = (-1)^c  \left|\left(
\begin{matrix}{i_1 \choose -\kappa_1} & \cdots & {i_d \choose
-\kappa_d} \cr {i_1 \choose 1-\kappa_1} & \cdots & {i_d \choose
1-\kappa_d} \cr \vdots & & \vdots \cr {i_1 \choose d-1-\kappa_1} &
\cdots & {i_d \choose d-1-\kappa_d}
\end{matrix}
\right)\cdot\left(
\begin{matrix}{n+1-j_d \choose -\gamma_d} & \cdots
& {n+1-j_1 \choose -\gamma_1} \cr {n+1-j_d \choose 1-\gamma_d} &
\cdots & {n+1-j_1 \choose 1-\gamma_1} \cr \vdots & & \vdots \cr
{n+1-j_d \choose d-1-\gamma_d} & \cdots & {n+1-j_1 \choose
d-1-\gamma_1}
\end{matrix}
\right)\right|,$$ where $\kappa_q := $ $\# \{ \t_p \mid \t_p > i_q
\}$, for $q = 1, \ldots, d$, and $\gamma_q := $ $\# \{ \t_p \mid
\t_p < j_q \}$, for $q = 1, \ldots, d$, and $c=\kappa_1 + \cdots +
\kappa_d + \gamma_1 + \cdots + \gamma_d$.
\end{thm}
\begin{proof} Follows immediately from Theorems \ref{prodmult},
\ref{ros-zel}, and Lemma \ref{oppreg}, in view of the fact that
$w_0\t=(n+1-\t_d,\ldots,n+1-\t_1)$ and
$w_0v=(n+1-j_d,\ldots,n+1-j_1)$.
\end{proof}

\end{document}